\documentclass{amsart}

\usepackage{psfrag,amssymb}

\newtheorem{theorem}{Theorem}
\newtheorem{lemma}[theorem]{Lemma}

\newtheorem{example}{Example}

\renewcommand{\le}{\leqslant}
\renewcommand{\ge}{\geqslant}
\renewcommand{\epsilon}{\varepsilon}
\def\({\left(}
\def\){\right)}
\DeclareMathOperator{\diam}{diam} 
\def\metric{\mathsf d}
\def\plus{\rtimes}

\begin{document}


\title{Markov partitions reflecting the geometry of~$\times2,\times3$}

\author{Thomas Ward}
\address{School of Mathematics, University of East Anglia, Norwich, CV4~7AL,
UK}
\author{Yuki Yayama}
\address{Centro de Modelamiento Matem{\'a}tico, Av. Blanco Encalada 2120 Piso
7, Santiago de Chile}

\thanks{This research was supported by EPSRC grant EP/C015754/1. The
authors thank Richard Miles for useful comments.}


\begin{abstract}
We give an explicit geometric description of the~$\times2,\times3$
system, and use this to study a uniform family of Markov partitions
related to those of Wilson and Abramov. The behaviour of these
partitions is stable across expansive cones and transitions in this
behaviour detects the non-expansive lines.
\end{abstract}

\maketitle

\section{Introduction}\label{sec:introduction}

Markov partitions are a powerful tool in the study of hyperbolic
diffeomorphisms of manifolds. Explicit constructions of Markov
partitions for hyperbolic toral automorphisms of
the~$2$-torus~$\mathbb T^2$ in the work of Adler and
Weiss~\cite{MR0257315} are an important paradigmatic example, and in
special situations the tight connection between the geometry of the
map and the partition found in~\cite{MR0257315} is extended to
automorphisms of~$\mathbb T^d$ with~$d>2$ by
Manning~\cite{MR1895206}. On the other hand, maps of objects that
are not quite manifolds arise naturally in dynamics, notably as
attractors of smooth maps in work of Bowen~\cite{MR0482842} and
Williams~\cite{MR0348794}. Thus a natural extension of the classical
theory of smooth maps of compact manifolds concerns maps of
solenoids; a useful overview and the history may be found in a paper
of Takens~\cite{MR2164385}. The simplest solenoids are algebraic:
compact groups that are locally isometric to products of local
fields.

The structure of a tangent space comprising a product of local
fields including non-Archimedean ones may be used to study various
dynamical properties of automorphisms of solenoids: exotic
orbit-growth properties by Chothi, Everest, Miles, Stevens and the
first author~\cite{MR1461206},~\cite{MR2339472}; entropy and
structure of~$\mathbb Z^d$-actions of entropy rank one by Einsiedler
and Lind~\cite{MR2031042}; topological entropy by Lind and the first
author~\cite{MR961739},~\cite{MR1882488}.

Our purpose here is to study geometrically natural Markov partitions
like those used by Abramov~\cite{MR0123346} and
Wilson~\cite{MR0390181} for one of the simplest examples in which
non-Archimedean directions arise in the tangent space, and to study
how the structure of those partitions changes in expansive cones.
This gives a simple geometrical instance of the `subdynamics
philosophy' of Boyle and Lind~\cite{MR1355295}. A combinatorial
instance of the same kind of structure appears in work of Miles and
the first author~\cite{MR2279271}, where it is shown that
directional zeta functions detect the non-expansive set for systems
of entropy rank one.

In order to do this, we describe the structure of the space obtained
by taking the invertible extension of the~$\mathbb N^2$-action
generated by~$x\mapsto2x\pmod{1}$ and~$x\mapsto3x\pmod{1}$ on the
additive circle in a geometric way. To simplify matters we
concentrate on this one example: the same kind of construction works
in those systems of entropy rank one with an adelic covering space,
but is significantly more involved. In principle the Markov and
generating properties of the partitions can be shown from our
geometric description, but for brevity we deduce some of these
properties from Wilson's results.

\section{The geometry of~$\times2,\times3$}

We make use of a simple version of the adelic machinery; an elegant
account may be found in Weil~\cite[Ch.~4]{MR0234930}. We wish to
describe the group~$X=\widehat{\mathbb Z[\frac{1}{6}]}$ of
characters on~$\mathbb Z[\frac{1}{6}]$ and its metric structure:
this group carries the automorphisms~$\alpha^{(1,0)}$
and~$\alpha^{(0,1)}$ dual to the automorphisms~$x\mapsto 2x$
and~$x\mapsto3x$ on~$\mathbb Z[\frac{1}{6}]$, and is a presentation
of the invertible extension of the~$\mathbb N^2$ action generated
by~$x\mapsto2x\pmod{1}$ and~$x\mapsto3x\pmod{1}$ on~$\mathbb T$. For
a prime~$p$, define the local field~$\mathbb Q_p$ to be the set of
formal power series~$\sum_{n\ge k}a_np^n$ with
digits~$a_n\in\{0,1,\dots,p-1\}$ and~$k\in\mathbb Z$, and with the
non-Archimedean metric~$\vert\cdot\vert_p$ induced by the~$p$-adic
norm~$\vert\sum_{n\ge k}a_np^n\vert_p=p^{-k}$ if~$a_k\neq0$. Notice
that~$\mathbb Q$ is a subfield of each~$\mathbb Q_p$ and
each~$\mathbb Q_p$ has a maximal compact subring~$\mathbb
Z_p=\{x\in\mathbb Q_p\mid\vert x\vert_P\le1\}$.

The homomorphism
\begin{eqnarray*}
\Delta:\mathbb Z[\textstyle\frac{1}{6}]&\longrightarrow&\mathbb
R\times\mathbb
Q_2\times\mathbb Q_3\\
r&\longmapsto&(r,r,r)
\end{eqnarray*}
embeds~$\mathbb Z[\frac{1}{6}]$ as a discrete (and hence closed)
subgroup of~$\mathbb R\times\mathbb Q_2\times\mathbb Q_3$ with
respect to the metric~$\metric(x,y)=\max\left\{\vert
x_{\infty}-y_{\infty}\vert,\vert x_2-y_2\vert_2,\vert
x_3-y_3\vert_3\right\}$, where we
write~$x=(x_{\infty},x_2,x_3)\in\mathbb R\times\mathbb
Q_2\times\mathbb Q_3$. Write~$G=\mathbb R\times\mathbb
Q_2\times\mathbb Q_3$ and~$\Gamma=\Delta(\mathbb Z[\frac{1}{6}])$.
The group~$X$ is the quotient~$G/\Gamma$ (this may be seen from
Weil~\cite[Ch.~4]{MR0234930}), and we wish to describe this quotient
space in a concrete way. In order to motivate this, notice that a
toral automorphism may be constructed as follows. The identity map
embeds~$\mathbb Z^d$ as a discrete subgroup of~$\mathbb R^d$, and a
choice of coset representatives for~$\mathbb R^d/\mathbb Z^d$ gives
an explicit geometric description of the map induced on the torus by
any automorphism of~$\mathbb R^d$ preserving~$\mathbb Z^d$. In order
to make this note self-contained and to rehearse the kind of
calculation needed later, we include the proof of the following two
lemmas, which are simple instance of a well-known principle (see
Weil~\cite[Ch.~4]{MR0234930} or Hewitt and Ross~\cite[\S~II.10,
Th.~10.15]{MR551496}).

\begin{lemma}
The set~$F=[0,1)\times\mathbb Z_2\times\mathbb Z_3$ is a fundamental
domain for~$\Gamma$ in~$G$.
\end{lemma}

\begin{proof}
The first step is to check that~$F$ is \emph{big} enough:
given~$x\in G$, can we find~$\gamma=(r,r,r)\in\Gamma$
with~$x-\gamma\in F$? To do this, write~$\{\sum_{n\ge
k}a_np^n\}=\sum_{n=k}^{-1}a_np^n$ for the fractional part
of~$x\in\mathbb Q_p$;~$\{t\}$ for the fractional part and~$\lfloor
t\rfloor$ for the integer part of~$t\in\mathbb R$. A calculation
shows that if
\[
r=\{x_2\}+\{x_3\}+\lfloor(x_{\infty}-\{x_2\}-\{x_3\}\rfloor
\]
then~$r\in\mathbb Z[\frac{1}{6}]$ and~$(x_{\infty}-r,x_2-r,x_3-r)\in
F$ as required.

The second step is to check that~$F$ is \emph{small} enough:
if~$x,y\in F$ define the same coset of~$\Gamma$ then they are equal.
Assume therefore that~$x,y\in F$ and~$x-y=(r,r,r)\in\Gamma$.
Then~$x_2-y_2\in\mathbb Z_2\cap\mathbb Z[\frac{1}{6}]=\mathbb
Z[\frac{1}{3}]$ and~$x_3-y_3\in\mathbb Z_3\cap\mathbb
Z[\frac{1}{6}]=\mathbb Z[\frac{1}{2}]$, so~$r\in\mathbb
Z[\frac{1}{3}]\cap\mathbb Z[\frac{1}{2}]=\mathbb Z$, and
therefore~$\{x_{\infty}\}=\{y_{\infty}\}$,
so~$x_{\infty}=y_{\infty}$ and~$r=0$ as required.
\end{proof}

This means that there is a bijection~$G/\Gamma\longleftrightarrow
F$; to go further we need to describe the image of the group
operation on~$G/\Gamma$ under this bijection.

\begin{lemma}\label{lemma:additionrule}
For~$s,t\in G$,
\[
(t+\Gamma)+(s+\Gamma)=\(\{t_{\infty}+s_{\infty}\},t_2+s_2-\lfloor
t_{\infty}+s_{\infty}\rfloor,t_3+s_3-\lfloor
t_{\infty}+s_{\infty}\rfloor\)+\Gamma
\]
is the unique coset representative for~$t+s$ in~$F$.
\end{lemma}

\begin{proof}
We wish to find the unique~$u\in F$ with the property that there is
some~$(r,r,r)$ in~$\Gamma$ with~$u=t+s-r$. We must
have~$u_{\infty}=\{t_{\infty}+s_{\infty}\}$, which forces~$r$ to
be~$\lfloor t_{\infty}+s_{\infty}\rfloor$; notice that we also then
have
\[
u_2=t_2+s_2-\lfloor t_{\infty}+s_{\infty}\rfloor\in\mathbb Z_2
\]
and
\[
u_3=t_3+s_3-\lfloor t_{\infty}+s_{\infty}\rfloor\in\mathbb Z_3
\]
since~$\mathbb Z_2,\mathbb Z_3$ are rings.
\end{proof}

Lemma~\ref{lemma:additionrule} may be written as follows: the
operation
\begin{equation}\label{eqn:additionformula}
t\plus s=\(\{t_{\infty}+s_{\infty}\},t_2+s_2-\lfloor
t_{\infty}+s_{\infty}\rfloor,t_3+s_3-\lfloor
t_{\infty}+s_{\infty}\rfloor\)
\end{equation}
makes~$F$ into a group~$X=(F,\plus)$ isomorphic to~$G/\Gamma$. An
explicit metric on~$X$ is given by
\[
\overline{\metric}(x+\Gamma,y+\Gamma)=\min_{r\in\mathbb
Z[\frac{1}{6}]}\max\{\vert x_{\infty}-y_{\infty}+r\vert_{\infty},
\vert x_{2}-y_{2}+r\vert_{2},\vert x_{3}-y_{3}+r\vert_{3}\}.
\]
Wilson~\cite{MR0390181} describes the same solenoid in a different
way, as a projective limit of circles
\begin{equation}\label{equation:wilsonform}
X\cong\{z\in\mathbb T^{\mathbb N_0}\mid 6z_{k+1}=z_k\pmod{1}\mbox{
for all }k\ge1\};
\end{equation}
points~$z,z'$ in this description are close if their
coordinates~$z_k,z'_k$ are close in~$\mathbb T$ for~$1\le k\le K$
for large~$K$. The isomorphism in~\eqref{equation:wilsonform} may be
thought of as follows. A given point~$z=(z_k)_{k\ge0}$ in the
right-hand side of~\eqref{equation:wilsonform} defines an
element~$z_0\in\mathbb T$; each choice of~$z_{k+1}$ given~$z_k$
defines a unique pair~$x_2^{(k)}\in\{0,1\}$
and~$x_3^{(k)}\in\{0,1,2\}$
satisfying~$z_{k+1}=\frac16z_{k}+\frac{x_2^{(k)}}{2}+\frac{x_3^{(k)}}{6}$
(thinking of~$z_{k+1}$ as a real number in~$[0,1)$). The isomorphism
is then defined by sending~$z$ to the
point~$\(z_0,\sum_{k\ge0}x_2^{(k)}2^k,\sum_{k\ge0}x_3^{(k)}3^k\)\in
X$. This isomorphism respects the metric structures (nearby points
in~$X$ correspond to nearby points in the projective limit) and is
equivariant with respect to the automorphisms we study. The
automorphisms~$\alpha^{(1,0)}:x\mapsto 2x$
and~$\alpha^{(0,1)}:x\mapsto3x$ on~$G$ preserve~$\Gamma$ and
therefore define automorphisms of~$X=(F,\plus)$.

To see how the group~$X$ works, we compute the
automorphisms~$\alpha^{(0,1)}$ (multiplication
by~$3$),~$\alpha^{(-1,0)}$ (multiplication by~$\frac{1}{2}$),
and~$\alpha^{(-1,1)}$ (multiplication by~$\frac{3}{2}$) explicitly.
By~\eqref{eqn:additionformula},
\[
\alpha^{(0,1)}(x)=x\plus x\plus x=\(\{3x_{\infty}\},3x_2-\lfloor
3x_{\infty}\rfloor,3x_3-\lfloor 3x_{\infty}\rfloor\).
\]
Notice that the map~$\alpha^{(0,1)}$ locally \emph{expands} the real
component by a factor of~$3$, locally \emph{contracts} the~$3$-adic
component by a factor of~$3$, and is a local isometry on
the~$2$-adic component.

Write~$x_p=\sum_{n\ge k}x_p^{(n)}p^n$ with
digits~$x_p^{(n)}\in\{0,1,\dots,p-1\}$ for~$n\ge k$. Then
\[
\alpha^{(-1,0)}(x)=\left(\textstyle\frac{1}{2}+\frac{1}{2}x_2^{(0)},\frac{1}{2}x_2+\frac{1}{2}x_2^{(0)},
\frac{1}{2}x_3+\frac{1}{2}x_2^{(0)}\right)
\]
(this is most easily verified by doubling the right-hand side).

Finally, by combining the two calculations we see
that~$\alpha^{(-1,1)}(x)$ is
\[
\(\left\{\textstyle\frac{3}{2}x_{\infty}+\frac{3}{2}x_2^{(0)}\right\},
\textstyle\frac{3}{2}x_2+\frac{3}{2}x_2^{(0)}-\left\lfloor\textstyle\frac{3}{2}x_{\infty}+\frac{3}{2}x_2^{(0)}
\right\rfloor,
\frac{3}{2}x_3+\frac{3}{2}x_2^{(0)}-\left\lfloor\textstyle\frac{3}{2}x_{\infty}+\frac{3}{2}x_2^{(0)}
\right\rfloor\).
\]
Locally the action of~$\alpha^{(a,b)}$ multiplies by~$2^a3^b$, and
therefore acts on each of the three directions in the tangent space
as shown in Table~\ref{table} ($u$, $s$ denote unstable and stable
directions).

\begin{table}
\begin{center}
\caption{\label{table}Stable and unstable directions.}
\begin{tabular}{|c|c|c|c|}
\hline
{region}&{$\mathbb R$}&$\mathbb Q_2$&$\mathbb Q_3$\\
\hline $a>0,b>0$&$u$&$s$&$s$\\
\hline $a<0,b>0,2^a3^b>1$&$u$&$u$&$s$\\
\hline $a>0,b<0,2^a3^b>1$&$u$&$s$&$u$\\
\hline $a<0,b<0$&$s$&$u$&$u$\\
\hline $a>0,b<0,2^a3^b<1$&$s$&$s$&$u$\\
\hline $a<0,b>0,2^a3^b<1$&$s$&$u$&$s$\\\hline
\end{tabular}
\end{center}
\end{table}

The first three regions shown in Table~\ref{table} are the expansive
regions in the sense of~\cite{MR1355295} and~\cite{MR1869066}
(expansive regions are defined in the Grassmannian space of lines
in~$\mathbb R^2$, of which the circle is a two-fold cover; the table
shows the six regions in the cover). There are three non-expansive
lines~$a=0$ (containing maps like~$\alpha^{(0,1)}$, which behaves
like an isometry on the~$2$-adic direction),~$b=0$ (containing maps
like~$\alpha^{(1,0)}$, which behaves like an isometry on
the~$3$-adic direction) and~$2^a3^b=1$ (which does not contain any
lattice points, but has a sequence of lattice points~$(a_k,b_k)$
converging to it with the property that the real Lyapunov
exponent~$\log\vert 2^{a_k}3^{b_k}\vert$ of the
map~$\alpha^{(a_k,b_k)}$ converges to zero as~$k\to\infty$).

\section{Stable Markov partitions}

It is clear that there cannot be a single finite partition that is
generating for all the maps~$\alpha^{(a,b)}$ as~$(a,b)$ varies
inside an expansive cone because the set of topological entropies of
the maps in a cone is unbounded. Thus, what we mean by ``stable'' is
that the Markov partition for~$\alpha^{(a,b)}$ is constructed in a
uniform manner across all~$(a,b)\in\mathbb Z^2$. We will see later
that the geometry of how the map acts on an atom of the partition is
uniform across each expansive cone but changes at each non-expansive
direction.

Recall that the na\"{\i}ve height (in the sense of Diophantine
geometry) of a non-zero rational~$r/s$ is defined to
be~$H(r/s)=\max\{\vert r\vert,\vert s\vert\}$. Thus Abramov's
formula~\cite{MR0123346} for the entropy of an automorphism of a
one-dimensional solenoid may be written~$h(T)=\log H(r/s)$ if~$T$ is the
map dual to multiplication by~$r/s$.

\begin{theorem}\label{theorem:stablemarkovpartitions}
For each~$(a,b)\in\mathbb Z^2\setminus\{(0,0)\}$
let~$\xi^{(a,b)}$ denote the partition
\[
\left\{
A_j=\left[\textstyle\frac{j}{H(2^a3^b)},\frac{j+1}{H(2^a3^b)}\right)\times\mathbb
Z_2\times\mathbb Z_3\mid 0\le j<H(2^a3^b) \right\}.
\]
Then~$\{\xi^{(a,b)}\}$ is a stable family of Markov partitions whose
geometry detects the non-expansive directions of~$\alpha$. The
partition~$\xi^{(a,b)}$ is generating for~$\alpha^{(a,b)}$ if and
only if~$\alpha^{(a,b)}$ is expansive.
\end{theorem}

The theory of Markov partitions in the topological (rather than
smooth) setting is developed by Adler~\cite{MR1477538}; by `Markov'
we mean that the partition obtained from~$\xi^{(a,b)}$ by using open
intervals in the real coordinate instead of half-open intervals
satisfies~\cite[Def.~6.1]{MR1477538}. Much of the proof of
Theorem~\ref{theorem:stablemarkovpartitions} will use results from
Wilson~\cite{MR0390181} that conceal the geometry of the actions. In
order to see how the maps act geometrically, we illustrate the
result by describing the partition and the action of the map on the
partition in some representative directions. In each figure the
image of the atom~$A_0$ of the partition is shaded.

\begin{example}\label{example:one}\rm
Consider the direction~$(1,0)$, with corresponding map
\[
\alpha^{(1,0)}(x)=\(\{2x_{\infty}\},2x_2-\lfloor
2x_{\infty}\rfloor,2x_2-\lfloor 2x_{\infty}\rfloor\).
\]
The partition~$\xi^{(1,0)}$ simply divides the real coordinate
into~$[0,\frac12)$ and~$[\frac12,1)$. We compute
\[
\alpha^{(1,0)}(\xi^{(1,0)})=
\{
[0,1)\times2\mathbb Z_2\times\mathbb Z_3,
[0,1)\times(1+2\mathbb Z_2)\times\mathbb Z_3
\}
\]
and
\[
\alpha^{(2,0)}(\xi^{(1,0)})=
\{
[0,1)\times(4\mathbb Z_2\cup1+4\mathbb Z_2)\times\mathbb Z_3,
[0,1)\times(2+4\mathbb Z_2\cup3+4\mathbb Z_2)\times\mathbb Z_3
\}.
\]
Similarly,
\[
\alpha^{(-1,0)}(\xi^{(1,0)})= \{
([0,\textstyle\frac14)\cup[\frac12,\frac34))\times\mathbb
Z_2\times\mathbb Z_3,
([\frac14,\frac12)\cup[\frac34,1))\times\mathbb Z_2\times\mathbb Z_3
\}
\]
and~$\alpha^{(-2,0)}(\xi^{(1,0)})$ is the partition into the sets
\[
\(
[0,\textstyle\frac18)\cup[\frac14,\frac38)\cup[\frac12,\frac58)\cup[\frac34,\frac78)
\)\times\mathbb Z_2\times\mathbb Z_3
\]
and
\[
\([\textstyle\frac18,\frac14)\cup[\frac38,\frac12)\cup[\frac58,\frac34)\cup[\frac78,1)
\)\times\mathbb Z_2\times\mathbb Z_3.
\]
\begin{figure}[htbp]
      \begin{center}
      \psfrag{a}{$\alpha^{(1,0)}$}
      \psfrag{b}{$\alpha^{(1,0)}$}
      \psfrag{r}{$[0,1)$}
      \psfrag{z2}{$\mathbb Z_2$}
      \psfrag{z3}{$\mathbb Z_3$}
      \psfrag{0}{$0$}
      \psfrag{1}{$1$}
      \psfrag{A0}{$A_0$}
\scalebox{1}{\includegraphics{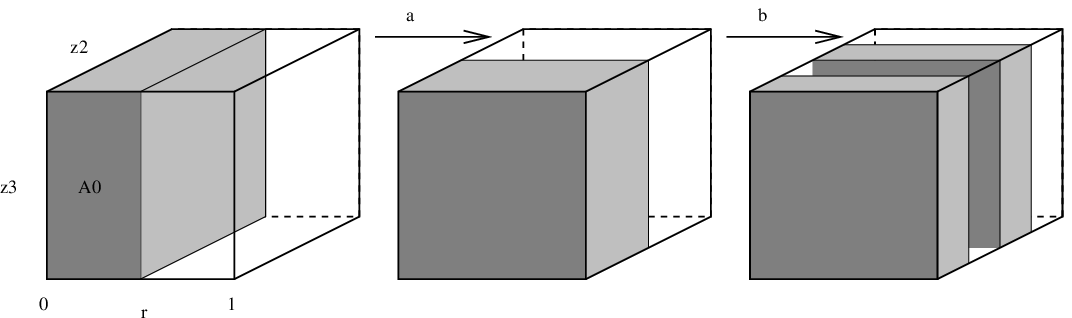}}
         \caption{\label{fig:multby2a} $\xi^{(1,0)}$, $\alpha^{(1,0)}(\xi^{(1,0)})$ and
         $\alpha^{(2,0)}(\xi^{(1,0)})$.}
\end{center}
\end{figure}

\begin{figure}[htbp]
      \begin{center}
      \psfrag{a}{$\alpha^{(-1,0)}$}
\scalebox{1}{\includegraphics{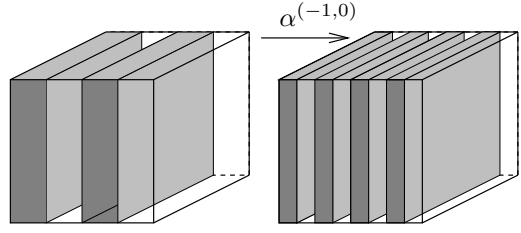}}
         \caption{\label{fig:multby2b} $\alpha^{(-1,0)}(\xi^{(1,0)})$ and
         $\alpha^{(-2,0)})(\xi^{(1,0)})$.}
\end{center}
\end{figure}
\noindent These partitions are illustrated in
Figure~\ref{fig:multby2a} for the forward direction and
Figure~\ref{fig:multby2b} for the reverse direction. Notice
that~$\bigvee_{n=-\infty}^{\infty}\alpha^{(n,0)}(\xi^{(1,0)})$ does
not separate the~$\mathbb Z_3$ coordinate, so the partition is not
generating for~$\alpha^{(1,0)}$. However, this does show that the
system~$(X,\alpha^{(1,0)})$ may be realized as an isometric
extension of a base system (which is an almost~$1$-$1$ image of a
full Bernoulli $2$-shift) by~$\mathbb Z_3$.
\end{example}

\begin{example}\rm
The expansive region~$ab>0$ is particularly simple because the
system~$(X,\alpha^{(a,b)})$ is (at each point with~$a>0,b>0$) simply
the invertible extension of the map~$x\mapsto2^a3^bx\pmod{1}$ on the
circle, and~$\xi^{(a,b)}$ is the usual partition into intervals of
width~$\frac{1}{2^a3^b}$ on~$[0,1)$ lifted to~$X$. The action
of~$\alpha^{(1,1)}$ (multiplication by~$6$) is illustrated in
Figure~\ref{fig:multby6a} for the forward direction and
Figure~\ref{fig:multby6b} for the reverse direction.
\begin{figure}[htbp]
      \begin{center}
      \psfrag{a}{$\alpha^{(1,1)}$}
      \psfrag{b}{$\alpha^{(1,1)}$}
\scalebox{1}{\includegraphics{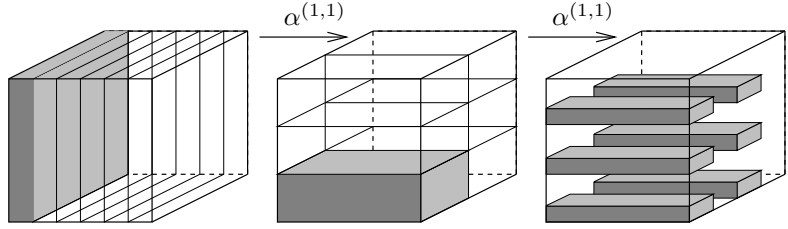}}
         \caption{\label{fig:multby6a} $\xi^{(1,1)}$, $\alpha^{(1,1)}(\xi^{(1,1)})$ and
         $\alpha^{(2,2)}(\xi^{(1,1)})$.}
\end{center}
\end{figure}

\begin{figure}[htbp]
      \begin{center}
      \psfrag{a}{$\alpha^{(-1,0)}$}
\scalebox{1}{\includegraphics{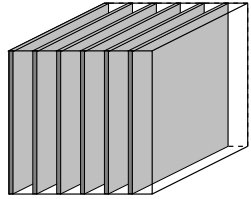}}
         \caption{\label{fig:multby6b} $\alpha^{(-1,-1)}(\xi^{(1,1)})$.}
\end{center}
\end{figure}
\end{example}

\begin{example}\rm
Now consider the map~$\alpha^{(-1,1)}$ (multiplication
by~$\frac32$). For this map the real and the~$2$-adic directions are
unstable and the~$3$-adic direction is stable. The
partition~$\xi^{(-1,1)}$ divides the real coordinate into three
pieces. A calculation shows that~$\alpha^{(-1,1)}(\xi^{(-1,1)})$
consists of the sets
\[
[0,\textstyle\frac12)\times\mathbb Z_2\times3\mathbb Z_3\cup
[\frac12,1)\times\mathbb Z_2\times(3\mathbb Z_3+2),
\]
\[
[\textstyle\frac12,1)\times\mathbb Z_2\times3\mathbb Z_3\cup
[0,\frac12)\times\mathbb Z_2\times(3\mathbb Z_3+1),
\]
and
\[
[0,\textstyle\frac12)\times\mathbb Z_2\times(3\mathbb Z_3+2)\cup
[\frac12,1)\times\mathbb Z_2\times(3\mathbb Z_3+1).
\]
The image of~$A_0$ under the maps~$\alpha^{(-1,1)}$
and~$\alpha^{(1,-1)}$ are shown in Figure~\ref{fig:multby3over2a}.

\begin{figure}[htbp]
      \begin{center}
      \psfrag{b}{$\alpha^{(-1,1)}$}
      \psfrag{a}{$\alpha^{(1,-1)}$}
\scalebox{1}{\includegraphics{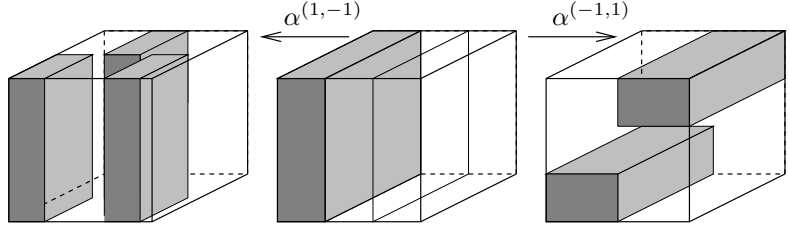}}
         \caption{\label{fig:multby3over2a} $\alpha^{(1,-1)}(\xi^{(-1,1)})$,
         $\xi^{(-1,1)}$ and $\alpha^{(-1,1)}(\xi^{(-1,1)})$.}
\end{center}
\end{figure}

\noindent A similar calculation shows
that~$\alpha^{(1,-1)}(\xi^{(-1,1)})$ consists of the sets
\[
[0,\textstyle\frac29)\times2\mathbb Z_2\times\mathbb Z_3\cup
[\frac13,\textstyle\frac59)\times(1+2\mathbb Z_2)\times\mathbb Z_3\cup
[\frac23,\textstyle\frac89)\times2\mathbb Z_2\times\mathbb Z_3,
\]
\[
[\textstyle\frac29,\textstyle\frac49)\times2\mathbb Z_2\times\mathbb Z_3\cup
[\frac59,\textstyle\frac79)\times(1+2\mathbb Z_2)\times\mathbb Z_3\cup
\(
[\frac89,\textstyle1)\times2\mathbb Z_2\times\mathbb Z_3\cup
[0,\textstyle\frac19)\times(1+2\mathbb Z_2)\times\mathbb Z_3
\),
\]
and the complement of their union. Notice that (for
example)~$\alpha^{-1}A_0\cap A_1\cap\alpha A_0$ does not consist of
a single rectangle.
\end{example}

\begin{proof}[{Proof of Theorem~\ref{theorem:stablemarkovpartitions} in the
region~$ab>0$.}]
Assume first that~$a>0$,~$b>0$, so that~$2^a3^b\in\mathbb N$, and
write~$\alpha=\alpha^{(a,b)}$,~$\xi=\xi^{(a,b)}$ throughout; the
partitions~$\alpha^{-1}(\xi),\xi, \alpha(\xi)$ are illustrated in
Figures~\ref{fig:multby6a} and~\ref{fig:multby6b} with the image and
pre-image of~$A_0$ shaded for the case~$(a,b)=(1,1)$. We claim that
the combinatorics of a full shift on~$6$ symbols suggested by
Figures~\ref{fig:multby6a} and~\ref{fig:multby6b} is indeed the
case. This (and other steps flagged below) may in principle be
extracted from Wilson's paper~\cite{MR0390181} but we prove it here
to show how the map works. We first need to check that an atom of
the form
\[
\alpha(A_{i_1})\cap\cdots\cap\alpha^{n}(A_{i_{n}}),
\]
for any choice of~$i_1,\dots,i_{n}\in\{0,1,\dots,2^a3^b-1\}$, is a
rectangle of the shape
\[
[0,1)\times (t_n+2^{an}\mathbb Z_2)\times(s_n+3^{bn}\mathbb Z_3)
\]
with an explicit description of~$t_n\in\{0,1,\dots 2^{an}-1\}$
and~$s_n\in\{0,1,\dots,3^{bn}-1\}$.
In order to do this, we need some notation
for the sets arising as the map is iterated.
The first iteration is straightforward, and
we can write
\[
\alpha(A_k)=[0,1)\times(2^a\mathbb Z_2-k)\times(3^b\mathbb Z_3-k)
\]
for~$0\le k\le 2^a3^b-1.$ The next iteration is more
complicated, because the image involves reduction modulo~$\Gamma$.
We compute
\begin{equation}\label{equation:suggestionE1fromyuki}
\alpha^2(A_k)=\bigsqcup_{\ell_1=0}^{2^a3^b-1} A_{k,\ell_1},
\end{equation}
where
\[
A_{k,\ell_1}= [0,1)\times\(2^a3^b(2^a\mathbb Z_2-k)-\ell_1\)\times
\(2^a3^b(3^b\mathbb Z_3-k)-\ell_1\)
\]
($\sqcup$ denoting a disjoint union). Continue, arriving at the
notation
\begin{equation}\label{equation:suggestionE2fromyuki}
\alpha^n(A_k)=\bigsqcup_{\ell_1=0}^{2^a3^b-1} \cdots
\bigsqcup_{\ell_{n-1}=0}^{2^a3^b-1}A_{k,\ell_1,\dots,\ell_{n-1}}
\end{equation}
for~$n\ge2$, in which each~$A_{k,\ell_1,\dots,\ell_{n-1}}$ is a set
of the form
\[
[0,1)\times \( 2^{an}\mathbb
Z_2-C(k,\ell_1,\dots,\ell_{n-2})-\ell_{n-1} \) \times \(
3^{bn}\mathbb Z_3-C(k,\ell_1,\dots,\ell_{n-2})-\ell_{n-1} \)
\]
where
\[
C(k,\ell_1,\dots,\ell_{n-2})= k(2^{a}3^{b})^{n-1}+\ell_1
(2^{a}3^{b})^{n-2}+\ell_2(2^{a}3^{b})^{n-3}+\cdots+\ell_{n-2}2^{a}3^{b}.
\]
Using this description, we claim that an atom
in~$\bigvee_{j=1}^{n}\alpha^j(\xi)$ can be written in the form
\begin{equation}\label{equation:nameofatomforwardpositivequadrant}
\alpha(A_{i_1})\cap\alpha^2(A_{i_2})\cap\cdots\cap\alpha^n(A_{i_n})=
A_{i_n,i_{n-1},\dots,i_1}
\end{equation}
for~$n\ge2$ and some~$0\le i_j<2^a3^b$,~$1\le j\le n$
where the right-hand side is
defined as above.

We prove the claim
in~\eqref{equation:nameofatomforwardpositivequadrant} by induction
on the length~$n$ starting with~$n=2$.
Clearly~$\alpha^2(A_{i_2})\supseteq A_{i_2,i_1}$ by definition. Now
\begin{eqnarray*}
A_{i_2,i_1}&=&[0,1)\times\(2^{2a}\mathbb Z_2-i_22^a3^b-i_1\)\times\(3^{2b}\mathbb Z_3-i_22^a3^b-i_1\)\\
&\subseteq&[0,1)\times\(2^a\mathbb Z_2-i_1\)\times\(3^b\mathbb
Z_3-i_1\)=\alpha(A_{i_1})
\end{eqnarray*}
since~$i_22^a3^b\mathbb Z_2\subseteq 2^a\mathbb Z_2$, and similarly
for the other terms, so~$\alpha(A_{i_1})\supseteq A_{i_2,i_1}$. Thus
\[
\alpha(A_{i_1})\cap\alpha^2(A_{i_2})\supseteq A_{i_2,i_1}.
\]
We now claim that~$\alpha(A_{i_1})\cap\alpha^2(A_{i_2})=
A_{i_2,i_1}$ by using~\eqref{equation:suggestionE1fromyuki} and
showing that
\[
A_{i_2,\ell}\cap \alpha(A_{i_1})\neq\emptyset
\]
for~$0\le\ell<2^a3^b$ implies that~$\ell=i_1$. To see this, note
first that if~$A_{i_2,\ell}\cap \alpha(A_{i_1})\neq\emptyset$
then~$A_{i_2,\ell}\subset\alpha(A_{i_1})$. Suppose that there is
some~$i_1'\neq i_1$, both in~$\{0,\dots,2^a3^b-1\}$,
with~$A_{i_2,i_1'}\cap\alpha(A_{i_1})\neq\emptyset$.
Then~$i_1-i_1'=2^ak_1$ and~$i_1-i_1'=3^bk_2$ for
some~$k_1,k_2\in\mathbb Z$, so (since~$2$ and~$3$ are
coprime),~$i_1\equiv i_1'\pmod{2^a3^b}$ and therefore~$i_1=i_1'$.

Now assume that~\eqref{equation:nameofatomforwardpositivequadrant}
holds for~$n\le k$. First notice that
\[
\alpha^{k+1}(A_{i_{k+1}})\supset A_{i_{k+1},i_k,\dots,i_1},
\]
and we claim that
\begin{equation}\label{equation:inductivesteppositivequadrantforward}
A_{i_k,\dots,i_1}\supset A_{i_{k+1},i_k,\dots,i_1}.
\end{equation}
Since
\begin{eqnarray*}
A_{i_k,\dots,i_1}&=&[0,1)\times\(
2^{ak}\mathbb Z_2-i_k(2^a3^b)^{k-1}-i_{k-1}(2^a3^b)^{k-2}-\cdots- i_1 \)\\
&&\qquad\qquad\qquad\times \( 3^{bk}\mathbb
Z_3-i_k(2^a3^b)^{k-1}-i_{k-1}(2^a3^b)^{k-2}-\cdots-i_1  \),
\end{eqnarray*}
\begin{eqnarray*}
A_{i_{k+1},i_k,\dots,i_1}&=&[0,1)\times\(2^{a(k+1)}\mathbb Z_2-
i_{k+1}(2^a3^b)^{k}-i_k(2^a3^b)^{k-1}-\cdots-i_1\)\\
&&\qquad\qquad\times \(3^{b(k+1)}\mathbb Z_3-
i_{k+1}(2^a3^b)^{k}-i_k(2^a3^b)^{k-1}-\cdots-i_1\),
\end{eqnarray*}
and
\begin{eqnarray*}
2^{a(k+1)}\mathbb Z_{2}&\subseteq& 2^{ak}\mathbb Z_{2}\\
3^{b(k+1)}\mathbb Z_{3}&\subseteq& 3^{bk}\mathbb Z_{3}\\
i_{k+1}(2^a3^b)^{k}\mathbb Z_2&\subseteq&2^{ak}\mathbb Z_2,\\
i_{k+1}(2^a3^b)^{k}\mathbb Z_3&\subseteq&3^{bk}\mathbb Z_3,
\end{eqnarray*}
we have~\eqref{equation:inductivesteppositivequadrantforward}, and
therefore
\begin{equation}\label{equation:onedirectionpositivequadrant}
A_{i_k,\dots,i_1}\cap\alpha^{k+1}(A_{i_{k+1}})\supseteq
A_{i_{k+1},\dots,i_1}.
\end{equation}
We now claim that there is equality
in~\eqref{equation:onedirectionpositivequadrant}. To see this,
assume that there is a choice
of~$\ell_1,\dots,\ell_k\in\{0,\dots,2^a3^b-1\}$
with~$A_{i_{k+1},\ell_k,\dots,\ell_1}\cap
A_{i_{k},i_{k-1},\dots,i_1}\neq\emptyset.$
By~\eqref{equation:inductivesteppositivequadrantforward}, and noting
that~$A_{i_{k+1},\ell_{k},\dots,\ell_{1}}\cap A_{i_{k},
i_{k-1},\dots,i_{1}}\neq \emptyset$ implies that~$A_{i_{k+1},
\ell_{k},\dots,\ell_{1}}$ is a subset of~$A_{i_k, i_{k-1},
\dots,i_1}$, it follows that
\[
2^{a(k+1)}\mathbb Z_2-i_{k+1}(2^a3^b)^{k}-i_k(2^a3^b)^{k-1}-\cdots-i_1
\]
and
\[
2^{a(k+1)}\mathbb
Z_2-i_{k+1}(2^a3^b)^{k}-\ell_k(2^a3^b)^{k-1}-\cdots-\ell_1
\]
are both subsets of
\[
2^{ak}\mathbb Z_2-i_k(2^a3^b)^{k-1}-\cdots-i_1,
\]
and similarly for the~$\mathbb Z_3$ component. Thus
\[
(i_k-\ell_k)(2^a3^b)^{k-1}+(i_{k-1}-\ell_{k-1})(2^a3^b)^{k-2}+\cdots+(i_1-\ell_1)\equiv0\pmod{(2^a3^b)^k}.
\]
Reducing this identity modulo~$2^a3^b$ shows that~$i_1=\ell_1$,
reducing modulo~$(2^a3^b)^2$ shows that~$i_2=\ell_2$, and so on.
Using~\eqref{equation:suggestionE2fromyuki}, it follows that there
is equality in~\eqref{equation:onedirectionpositivequadrant} as
required,
proving~\eqref{equation:nameofatomforwardpositivequadrant}.

Now we consider an atom of the form
\[
A_{i_0}\cap\alpha^{-1}(A_{i_1})\cap\cdots\cap\alpha^{-n}(A_{i_n});
\]
we wish to prove a statement
like~\eqref{equation:nameofatomforwardpositivequadrant} for these
atoms, by showing that each such atom is a
rectangle of the form~$J\times\mathbb Z_2\times\mathbb Z_3$
for an explicitly described interval~$J$ of width~$\frac{1}{(2^a3^b)^{n+1}}$.
A calculation shows that
\[
\alpha^{-1}(A_k)=\bigsqcup_{\ell=0}^{2^a3^b-1}
\left[\textstyle\frac{k}{(2^a3^b)^2}+\textstyle\frac{\ell}{2^a3^b},
\frac{k+1}{(2^a3^b)^2}+\frac{\ell}{2^a3^b}\right)
\times\mathbb Z_2\times\mathbb Z_3
=\bigsqcup_{\ell=0}^{2^a3^b-1}A^{k,\ell},
\]
and in general we have
\begin{equation}\label{equation:suggestionfromyukiE3}
\alpha^{-n}(A_k)=\bigsqcup_{\ell_1=0}^{2^a3^b-1}\cdots\bigsqcup_{\ell_n=0}^{2^a3^b-1}
A^{k,\ell_1,\dots,\ell_n}
\end{equation}
for~$n\ge1$, with
\[
A^{k,\ell_1,\dots,\ell_n}= \left[
\textstyle\frac{k}{(2^a3^b)^{n+1}}+D(\ell_1,\dots,\ell_n),
\textstyle\frac{k+1}{(2^a3^b)^{n+1}}+D(\ell_1,\dots,\ell_n)\right)\times\mathbb Z_2\times\mathbb Z_3
\]
where
\[
D(\ell_1,\dots,\ell_n)=\frac{\ell_1}{(2^a3^b)^n}+\frac{\ell_2}{(2^a3^b)^{n-1}}+
\cdots+\frac{\ell_n}{2^a3^b}.
\]
We claim that
\begin{equation}\label{equation:reversepositivequadrantclaim}
A_{i_0}\cap\alpha^{-1}(A_{i_1})\cap\cdots\cap\alpha^{-n}(A_{i_n})=
A^{i_n,i_{n-1},\dots,i_0}.
\end{equation}
for~$n\ge1$. For~$n=1$,
\[
A_{i_0}\cap\alpha^{-1}(A_{i_1})\supseteq\left[\textstyle\frac{i_0}{2^a3^b},\frac{i_0+1}{2^a3^b}\right)
\times\mathbb Z_2\times\mathbb Z_3\cap A^{i_1,i_0}=A^{i_1,i_0}
\]
since~$[\frac{i_0}{2^a3^b},\frac{i_0+1}{2^a3^b})\supseteq
[\frac{i_1}{(2^a3^b)^2}+\frac{i_0}{2^a3^b},\frac{i_1+1}{(2^a3^b)^2}+\frac{i_0}{2^a3^b})$.
Thus~$A_{i_0}\cap\alpha^{-1}(A_{i_1})=A^{i_1,i_0}$ since the width
of the real interval defining~$A_{i_0}$ is~$\frac{1}{2^a3^b}$ and
by~\eqref{equation:suggestionfromyukiE3} the real coordinates of the
sets in~$\alpha^{-1}(A_{i_1})$ are intervals, each of
width~$\frac{1}{(2^a3^b)^2}$ and with the property that the left
end-points of distinct intervals are at least~$\frac{1}{2^a3^b}$
apart.

Now assume that~\eqref{equation:reversepositivequadrantclaim} holds
for~$n\le k$, so that~$\bigcap_{j=0}^{k+1}\alpha^{-j}(A_{i_j})$ can
be written as the intersection of
\[
\left[
D(i_k,\dots,i_0),\textstyle\frac{1}{(2^a3^b)^{k+1}}+D(i_k,\dots,i_0)
\right)\times\mathbb Z_2\times\mathbb Z_3=A^{i_{k},\dots,i_{0}}
\]
with
\[
\bigsqcup_{0\le
j_1,\dots,j_{k+1}<2^a3^b}A^{i_{k+1},j_1,\dots,j_{k+1}}.
\]
It follows that
\[
\bigcap_{j=0}^{k+1}\alpha^{-j}(A_{i_j})\negthinspace\supseteq\negthinspace\negthinspace
\left[D(i_{k+1},\dots,i_0),\textstyle\frac{1}{(2^a3^b)^{k+2}}+
D(i_{k+1},\dots,i_0)\right)\negthinspace\times\mathbb
Z_2\times\mathbb Z_3\negthinspace\negthinspace=\negthinspace
A^{i_{k+1},\dots,i_{0}}.
\]
Notice that the width of the real interval defining the
set~$A^{i_{k},\dots,i_{0}}$ is~$\frac{1}{(2^a3^b)^{k+1}}$. Now
by~\eqref{equation:suggestionfromyukiE3} each member of the real
projection of~$\alpha^{-(k+1)}(A_{i_{k+1}})$ has
length~$\frac{1}{(2^a3^b)^{k+2}}$ and each of these intervals has
the property that the left end-points of distinct intervals are at
least distance~$\frac{1}{(2^a3^b)^{k+1}}$ apart,
showing~\eqref{equation:reversepositivequadrantclaim} for~$n=k+1$
and hence for all~$n$ by induction.

By~\eqref{equation:nameofatomforwardpositivequadrant}
and~\eqref{equation:reversepositivequadrantclaim}, the atom
\[
\bigcap_{j=-n}^{n}\alpha^{j}(A_{i_j})= A^{i_{n},\dots,i_0}\cap
A_{i_n,\dots,i_1}
\]
is a single rectangle with real width~$\frac{1}{(2^a3^b)^{n+1}}$,
$2$-adic width~$\frac{1}{(2^a)^n}$ and~$3$-adic
width~$\frac{1}{(3^b)^n}$. It follows that~$\xi$ satisfies a strong
form of the condition~\cite[Exercise~6.1]{MR1477538}. Moreover,
\[
\diam\(\vphantom{\sum}\right.\negmedspace\bigvee_{j=-n}^{n}\alpha^{j}(\xi)\negmedspace
\left.\vphantom{\sum}\)\rightarrow0
\]
as~$n\to\infty$, so~$\xi$ is a generating Markov partition in the
sense of~\cite{MR1477538}.
\end{proof}

\begin{proof}[{Proof of Theorem~\ref{theorem:stablemarkovpartitions} in other
regions.}]
Away from the positive and negative quadrants~$ab>0$ the behaviour
of~$\xi=\xi^{(a,b)}$ under the map~$\alpha=\alpha^{(a,b)}$ is more
complicated. In particular, as seen in
Figure~\ref{fig:multby3over2a}, an atom in~$\xi\vee\alpha\xi$ need
not be a rectangle even in expansive directions. However, in an
expansive direction the partition~$\xi$ corresponds under the map
described after~\eqref{equation:wilsonform} to the
partition~$\pi_0^{-1}S(H(2^a3^b))$ used by
Wilson~\cite[Th.~2.4]{MR0390181}. Notice that for any~$(a,b)$ in an
expansive region, the group~$\Sigma_{mn}$ in the notation
of~\cite{MR0390181}, where~$\frac{m}{n}=2^a3^b$, is~$X$. Wilson
shows that this partition is a Bernoulli generator, so
\[
\bigcap_{j=0}^{n}\alpha^{j}(A_{i_j})\neq\emptyset,
\bigcap_{j=-n}^{0}\alpha^{j}(A_{i_j})\neq\emptyset\implies
\bigcap_{j=-n}^{n}\alpha^{j}(A_{i_j})\neq\emptyset
\]
(as in~\cite[Exercise~6.1]{MR1477538}); he also shows that an atom
in~$\bigvee_{j=-n}^{n}\alpha^{j}(\xi)$ lies inside a cylinder
defined by small intervals in many coordinates in the
description~\eqref{equation:wilsonform}, so
\[
\diam\( \bigvee_{j=-n}^{n}\alpha^{j}(\xi) \)\rightarrow0.
\]
It follows that~$\xi$ is a generating Markov partition
for~$\alpha^{(a,b)}$.

There are three non-expansive directions, but only two of them
contain non-trivial lattice points: Example~\ref{example:one} shows
that~$\xi^{(1,0)}$ is not generating under~$\alpha^{(1,0)}$; the
other direction~$(0,1)$ is similar.
\end{proof}

An impression of the complexity of a generating Markov partition may
be gained by comparing the dynamical zeta function of the resulting
symbolic cover shift map~$\sigma^{(a,b)}$ to the zeta function of
the original map~$\alpha^{(a,b)}$. In the positive
quadrant~$a>0,b>0$, where we have seen that the partition~$\xi$
behaves very simply, we
have~$\zeta_{\sigma^{(a,b)}}(z)=\frac{1}{1-H(a,b)z}$
while~$\zeta_{\alpha^{(a,b)}}(z)=\frac{1-z}{1-H(a,b)z}$ since only
one pair of points of each period are identified by the factor map
defined by the partition. In contrast, in the
region~$a<0,b>0,2^a3^b>1$ (for example) we
have~$\zeta_{\sigma^{(a,b)}}(z)=\frac{1}{1-3^bz}$
while~$\zeta_{\alpha^{(a,b)}}(z)=\frac{1-2^az}{1-3^bz}$, reflecting
the fact that more periodic points in the full~$3^b$-shift are
identified under the factor map. Finally, in a non-expansive
direction (like~$a=1,b=0$) the zeta function of~$\alpha^{(a,b)}$ is
not even a rational function (it is shown in~\cite{MR2180241} that
the zeta function has a natural boundary on the circle~$\vert z
\vert=\frac{1}{2}$ in this case; the influence on the zeta function
of further directions in which an automorphism of a solenoid acts
like an isometry is studied by Miles~\cite{MR2308145} and the first
author~\cite{MR1458718}).

\begin{figure}[htbp]
      \begin{center}
      \psfrag{Q3}{$\mathbb Q_3$}\psfrag{Q2}{$\mathbb
      Q_2$}\psfrag{R}{$\mathbb
      R$}\psfrag{line}{$2^x3^y=1$}\psfrag{x2}{$\times
      2$}\psfrag{x3}{$\times 3$}
\scalebox{1}{\includegraphics{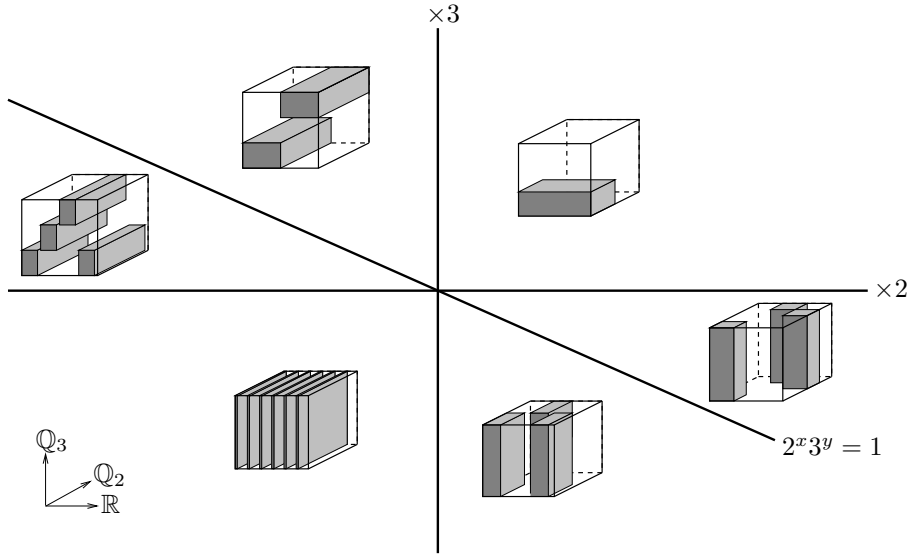}}
         \caption{\label{fig:seven} Geometry of~$\alpha^{(a,b)}(A_0)$ in expansive cones.}
\end{center}
\end{figure}

The variation in geometrical properties of the
partition~$\xi^{(a,b)}$ across each expansive cone is illustrated in
Figure~\ref{fig:seven}: a representative shape
of~$\alpha^{(a,b)}(A_0)$ is shown shaded in each expansive cone. The
transitions across the axes are clear; at the line~$2^x3^y=1$ all
that changes is the sign of the real Lyapunov exponent.


\providecommand{\bysame}{\leavevmode\hbox
to3em{\hrulefill}\thinspace}
\providecommand{\MR}{\relax\ifhmode\unskip\space\fi MR }
\providecommand{\MRhref}[2]{%
  \href{http://www.ams.org/mathscinet-getitem?mr=#1}{#2}
} \providecommand{\href}[2]{#2}

\end{document}